\documentclass[a4paper,12pt]{article}

\usepackage[latin1]{inputenc}
\usepackage{amsmath}
\usepackage{amsfonts}
\usepackage{amssymb}
\usepackage{amsthm}
\usepackage{fontenc}
\usepackage[pdftex]{graphicx}
\usepackage[left=2.5cm,right=2.5cm,top=2cm,bottom=2cm]{geometry}
\usepackage{cite}

\newtheorem{theorem}{Theorem}
\newtheorem{lemma}{Lemma}

\newcommand{\C}{\mathbb{C}}

\newcommand{\R}{\mathbb{R}}

\newtheorem{corollary}{Corollary}

\newcommand{\ux}{\underline{x}}
\newcommand{\uy}{\underline{y}}

\begin{document}

\title{On Riemann problems for monogenic functions in lower dimensional non-commutative Clifford algebras}

\author{Carlos Daniel Tamayo-Castro$^{1}$, Ricardo Abreu-Blaya$^{2}$\\ and\\ Juan Bory-Reyes$^{3}$}
\date {\small{$^{1}$ Instituto de Matem\'aticas. Universidad Nacional Aut\'onoma de M\'exico, Ciudad M\'exico, M\'exico.\\
e-mail: cdtamayoc@comunidad.unam.mx\\
$^{2}$ Facultad de Matem\'atica, Universidad Aut\'onoma de Guerrero, Chilpancingo, Guerrero, M\'exico.\\
e-mail: rabreublaya@yahoo.es\\
$^{3}$ SEPI-ESIME-Zac. Instituto Polit\'ecnico Nacional, Ciudad M\'exico, M\'exico.\\
e-mail: juanboryreyes@yahoo.com}}

\maketitle

\begin{abstract} 
In this paper, we mainly consider the Riemann boundary value problems for lower dimensional non-commutative Clifford algebras valued monogenic functions. The solutions are given in an explicit way and concrete examples are presented to illustrate the results.

\vspace{0.3cm}

\small{
}  
\end{abstract}
\section{Introduction}
The theory of Riemann boundary value problem (RBVP for short) for analytic functions of one complex variable was extensively studied for many researchers so far, see for example \cite{Gajov, Lu} for extensive treatments and discussions. On the other hand, a very detailed treatment of the RBVP for generalized analytic functions, as well as for many other linear and nonlinear elliptic systems in the plane do appear in many sources, see e.g., \cite{B, GB, Ve, We} or elsewhere.

Besides the theoretical significance of the RBVP, their study is closely connected with the theory of singular integral equations \cite{Mu} and has a wide range of applications in almost all areas of physics and engineering such as for instance in electromagnetism, optic, elasticity, fluid dynamics, geophysics, theory of orthogonal polynomials, in asymptotic analysis up to modern quantum field theory.

From both, mathematical convenience and physical relevancy, in the Clifford algebras-based approach, there has been a great deal of activity in the generalization of the classical-posed Riemann boundary value problems to higher dimensions. Several interesting results about RBVP for monogenic functions and harmonic functions in Clifford analysis (or one of its current research frameworks) are presented in \cite{BTA, BADK, Br, GZ, LJ, ZG, LZ1, LZ2, LJD, L1, L2, MYKC, ShV, CKM} and the references therein. 

Our purpose is to describe a sufficiently complete picture of solvability of the RBVP for monogenic functions in lower dimensional non-commutative Clifford algebras mainly concern the case of the so-called paravectorial and vectorial Clifford analysis.

\section{Preliminaries and Notations}
In this section we recall necessary basic facts about Clifford analysis which will be needed in the sequel. The best standard reference here is \cite{BDS}.

The real Clifford algebra associated with $\mathbb{R}^{n}$ endowed with the Euclidean metric is the minimal enlargement of $\mathbb{R}^{n}$ to a real linear associative algebra $\mathbb{R}_{0, n}$ with identity such that $x^{2} = -|x|^{2}$ for any $x \in \mathbb{R}^{n}$. 

It thus follows that if $\{e_{j}\}_{j = 1}^{n}$ is the standard orthonormal basis of $\mathbb{R}^{n}$ then we must have that $e_{i}e_{j} + e_{j}e_{i} = -2\delta_{ij}$, with $\delta_{ij}$ the Kronecker delta. So, one denotes an arbitrary $a \in \mathbb{R}_{0, n}$ by $a = \sum_{A \subseteq N}a_{A}e_{A}$, $N = \{1, \ldots, n \}$, $a_{A} \in \mathbb{R}$ where $e_{\emptyset} = e_{0} = 1$, $e_{\{ j\}} = e_{j}$ and $e_{A} = e_{\beta_{1}}\cdots e_{\beta_{k}}$ for $A = \{\beta_{1}, \ldots, \beta_{k} \}$ where $\beta_{j} \in \{1, \ldots, n \}$ and $\beta_{1} < \ldots < \beta_{k}$. 

Conjugation in $\R_{0,n}$ is defined as the anti-involution $a\mapsto\overline{a}:= \sum_{A}a_{A}\overline{e}_{A}$ for which 
\begin{equation*}
\overline{e}_{A} := (-1)^{k}e_{\beta_{k}} \cdots e_{\beta_{2}}e_{\beta_{1}}.
\end{equation*}

Put $\mathbb{R}^{(k)}_{0, n} = span_{\mathbb{R}}(e_{A}: |A| = k)$. Then clearly $\mathbb{R}^{(k)}_{0, n}$ is a subspace of $\mathbb{R}_{0, n}$ (the $k$-vectors in this class) and 
\begin{equation*}
\mathbb{R}_{0, n} = \bigoplus_{k = 0}^{n}\mathbb{R}^{(k)}_{0, n}.
\end{equation*}
The projection operator of $\mathbb{R}_{0, n}$ on $\mathbb{R}^{(k)}_{0, n}$ is denoted by $[\bullet]_{k}$ and $\mathbb{R}$ and $\mathbb{R}^{n}$ will be identified with $\mathbb{R}^{(0)}_{0, n}$ and $\mathbb{R}^{(1)}_{0, n}$ respectively.

Let us highlight the important fact that $\mathbb{R}_{0, n} = \mathbb{R}^{+}_{0, n}\oplus e_{1}\mathbb{R}^{+}_{0, n}$, where 
\[\mathbb{R}^{+}_{0, n}:=\bigoplus_{k-even}\mathbb{R}^{(k)}_{0, n}.\]
Then if $a \in \mathbb{R}_{0, n}$ we have the decomposition 
\begin{equation}\label{eo}
a = a_{0} + e_{1}a_{1},
\end{equation}
where $a_{0}, a_{1}$ will be referred to as its even and odd part respectively.

An important subspace of the real Clifford algebra $\mathbb{R}_{0, n}$ is the so-called space of paravectors $\mathbb{R}^{(0)}_{0, n}\oplus\mathbb{R}^{(1)}_{0, n}$, being the sum of scalars and
vectors. Notice that for each $x = (x^{0}, x^{1}, \ldots, x^{n}) \in \mathbb{R}^{n + 1}$, to be identified with
\begin{equation*}
x = x^{0} + \sum_{j = 1}^{n}x^{j}e_{j} \in \mathbb{R}^{(0)}_{0, n}\oplus\mathbb{R}^{(1)}_{0, n},
\end{equation*}
there should hold that
\begin{equation*}
x\overline{x} = \overline{x}x = |x|^{2}.
\end{equation*}
For further use we shall be considering the one to one mapping
\begin{equation*}
\begin{array}{cccc}
\alpha: & \mathbb{R}\oplus\mathbb{R}_{0,n - 1}^{(1)}   & \rightarrow & \mathbb{R}_{0,n}^{(1)}    \\ 
& x^{1} + \sum_{i = 1}^{n - 1}x^{i + 1}e_{i} & \rightarrow  & \sum_{i = 1}^{n}x^{i}e_{i},   
\end{array} 
\end{equation*}
As well the isomorphism
\begin{equation*}
\begin{array}{cccc}
\beta: & \mathbb{R}^{+}_{0, n}   & \rightarrow & \mathbb{R}_{0, n - 1} \\
       &     e_{1}e_{i + 1}      & \rightarrow &  e_{i}.
\end{array} 
\end{equation*}

Classical Clifford analysis consists in set up a function theory defined on Euclidean space and take values in a real Clifford algebra. The function theory concentrates on the notion of monogenic functions belonging to the kernel of a generalized Cauchy-Riemann operator (paravectorial Clifford analysis), or to that of its vectorial part, that is, the Dirac operator (vectorial Clifford analysis). In this way, Clifford analysis may be considered both as a generalization to higher dimension of the theory of holomorphic functions in the complex plane and a refinement of classical harmonic analysis due to fact that these differential operators factorize the Laplacian. 

We start with the paravector Clifford analysis case. The considered functions $u$ are defined in $\Omega\subset\mathbb{R}^{n + 1}$ and take values in (subspace of) the real Clifford algebra $\mathbb{R}_{0, n}$. These functions may be written as
\begin{equation*}
u(x) = \sum_{A}u_{A}(x)e_{A},
\end{equation*}
where $u_{A}$ are $\mathbb{R}$-valued functions.
 
We say that $u$ belongs to some classical class of functions on $\Omega$  if each of its components $u_{A}$ belongs to that class.

The theory of paravectorial monogenic functions with values in Clifford algebras generalizes in a natural way the theory of holomorphic functions of one complex variable to the $(n + 1)$-dimensional Euclidean space. Monogenic functions are null solutions of the generalized Cauchy-Riemann operator in $\mathbb{R}^{n + 1}$ defined by
\begin{equation*}
D_{n} := \sum^{n}_{j = 0}e_{j}\dfrac{\partial}{\partial x_{j}}.
\end{equation*}
The left (right) fundamental solution of this first order elliptic operator is given by 
\begin{equation*}
\begin{array}{cc}
E_{n}(x) = \dfrac{1}{\sigma_{n + 1}}\dfrac{\overline{x}}{|x|^{n + 1}}, &  x\in\mathbb{R}^{n + 1}\setminus \{0\}
\end{array} 
\end{equation*}
where $\sigma_{n + 1}$ is the area of the unit sphere in $\mathbb{R}^{n + 1}$. Clearly, $E_{n}$ is both left and right monogenic in $\mathbb{R}^{n + 1}\setminus \{0\}$.

For $\Omega\subset\mathbb{R}^{n + 1}$ open and $u \in C^{1}(\Omega)$ then $u$ will be called left (resp. right) monogenic in $\Omega$ if $D_{n}u = 0$ (resp. $uD_{n} = 0$) in $\Omega$. Furthermore, for a non-open $\Omega$ we will call $u$ monogenic in $\Omega$ if it is monogenic in some open neighborhood of $\Omega$.

A powerful tool in obtaining basic examples of left monogenic functions is the left Cliffordian Cauchy type integral. Let $\Omega\subset\mathbb{R}^{n + 1}$ a simply connected bounded domain with a smooth boundary $L$. We let $d\nu$ denote the surface measure on $L$.

For each $u\in C(L)$ its left Cliffordian Cauchy type integral is defined by 
\begin{equation*}
\begin{array}{cc}
(\mathcal{C}_{L}u)(x) := \int_{L}E_{n}(y - x)\kappa(y)u(y)d\nu(y), & x\notin L,
\end{array} 
\end{equation*}
and its singular version, the singular Cauchy type integral (also called the Hilbert transform) on $L$ to be
\begin{equation*}
\begin{array}{cc}
(\mathcal{S}_{L}u)(x) := 2\int_{L}E_{n}(y - x)\kappa(y)(u(y) - u(x))d\nu(y) + u(x), & x\in L. 
\end{array} 
\end{equation*}
Hereby $\kappa(y)$ denotes the outward pointing unit normal to $L$ at $y \in L$ and the integral in $\mathcal{S}_{L}$ is evaluated in the sense of Cauchy principal value.

On the other hand basic examples of right monogenic functions are obtained by means of the right Cliffordian Cauchy type integral
\begin{equation*}
\begin{array}{cc}
(u\mathcal{C}_{L})(x) := \int_{L}u(y)\kappa(y)E_{n}(y - x)d\nu(y), & x\notin L.
\end{array} 
\end{equation*}
We can now led to the vectorial Clifford analysis situation. Here, $\ux = (x^{1}, \ldots, x^{n}) \in \mathbb{R}^{n}$ will be identified with
\begin{equation*}
\ux = \sum_{j = 1}^{n}x^{j}e_{j} \in \mathbb{R}^{(1)}_{0, n}.
\end{equation*}
Monogenic functions in this context means solutions of the Dirac operator in $\mathbb{R}^{n}$
\begin{equation*}
\partial_{n} := \sum^{n}_{j = 1}e_{j}\dfrac{\partial}{\partial x^{j}}.
\end{equation*}
The fundamental solution of the Dirac operator is given by
\begin{equation*}
\begin{array}{cc}
\vartheta_{n}(\ux) = \dfrac{1}{\sigma_{n}}\dfrac{\overline{\ux}}{|\ux|^{n}}, &  \ux\in\mathbb{R}^{n}\setminus \{0\}
\end{array} 
\end{equation*}
where $\sigma_{n}$ is the area of the unit sphere in $\mathbb{R}^{n}$. If $\Omega$ is open in $\mathbb{R}^{n}$ and $u \in C^{1}(\Omega)$ then $u$ is said to be left (resp. right) monogenic in $\Omega$ if $\partial_{n}u = 0$ (resp. $u\partial_{n} = 0$) in $\Omega$.

The corresponding Cliffordian Cauchy type integrals in vectorial Clifford analysis setting occur analogously with $\vartheta_{n}(\ux)$ in place of $E_{n}(x)$.

In what follows, for a given function
\begin{equation*}
\begin{array}{cccc}
u: & \mathbb{R}_{0,n}^{(1)}   & \rightarrow & \mathbb{R}^{+}_{0, n},   \\ 
\end{array} 
\end{equation*}
we define
\begin{equation*}
\begin{array}{cccc}
\widehat{u} : & \mathbb{R}\oplus\mathbb{R}_{0,n - 1}^{(1)}   & \rightarrow & \mathbb{R}_{0, n - 1},   \\ 
\end{array} 
\end{equation*}
by $\widehat{u}(x) := \beta \circ u \circ \alpha (\ux) = \beta(u(\alpha(\ux)))$, and $x = \alpha(\ux)$.
\section{Reduction procedure for the RBVP}
From now on we regard $\Omega$ as a simply connected and bounded domain of $\R^{n}$ with smooth boundary $L$. We use the temporary notation $\Omega_+:=\Omega$, $\Omega_-:=\R^n\setminus\{\Omega\cup L\}$.
 
The Riemann boundary value problem for monogenic functions in vectorial Clifford analysis may be described as follows:

Let two $\mathbb{R}_{0, n}-$valued functions $G, g$ belonging to $C^{0,\nu}(L)$. Find a function $\Phi$ monogenic on $\R^n\setminus L$ continuously extendable from $\Omega_{\pm}$ to $L$ such that the following condition of their boundary values $\Phi^{\pm}$ on $L$ holds
\begin{equation}\label{ProbOrig}
\Phi^{+}(\ux) - G(\ux)\Phi^{-}(\ux) = g(\ux),\,\ux\in L.
\end{equation}
with $G(\ux) \neq 0$.

After using the decompostion \eqref{eo}, we have
\begin{equation}\label{Decomp}
\begin{array}{c}
\Phi^{+}(\ux) = \Phi_{0}^{+}(\ux) + e_{1}\Phi_{1}^{+}(\ux)\\
\Phi^{-}(\ux) = \Phi_{0}^{-}(\ux) + e_{1}\Phi_{1}^{-}(\ux)\\
G(\ux) = G_{0}(\ux) + e_{1}\overline{}G_{1}(\ux)\\
g(\ux) = g_{0}(\ux) + e_{1}g_{1}(\ux).
\end{array}
\end{equation}
Substituting (\ref{Decomp}) into (\ref{ProbOrig}) yields
\begin{equation*}
[\Phi_{0}^{+}(\ux) -(G_{0}(\ux)\Phi_{0}^{-}(\ux) - G^{*}_{1}(\ux)\Phi_{1}^{-}(\ux))] + e_{1}[\Phi_{1}^{+}(\ux) - (G^{*}_{0}(\ux)\Phi_{1}^{-}(\ux) + G_{1}(\ux)\Phi_{0}^{-}(\ux))] = 
\end{equation*}
\begin{equation*}
= g_{0}(\ux) + e_{1}g_{1}(\ux),
\end{equation*}
where $G^{*}_{j}(\ux) = -e_{1}G_{j}(\ux)e_{1}, j = 0, 1$. 

So that we arrive to the system
\begin{equation}\label{System}
\Big\lbrace \begin{array}{c}
\Phi_{0}^{+}(\ux) -(G_{0}(\ux)\Phi_{0}^{-}(\ux) - G^{*}_{1}(\ux)\Phi_{1}^{-}(\ux)) = g_{0}(\ux), \\
\Phi_{1}^{+}(\ux) - (G_{1}(\ux)\Phi_{0}^{-}(\ux) + G^{*}_{0}(\ux)\Phi_{1}^{-}(\ux))= g_{1}(\ux).
\end{array}
\end{equation}
The system (\ref{System}) becomes 
\begin{equation*}
\Big\lbrace \begin{array}{c}
\widehat{\Phi}_{0}^{+}(\ux) -(\widehat{G}_{0}(\ux)\widehat{\Phi}_{0}^{-}(\ux) - \widehat{G}^{*}_{1}(\ux)\widehat{\Phi}_{1}^{-}(\ux)) = \widehat{g}_{0}(\ux), \\
\widehat{\Phi}_{1}^{+}(\ux) - (\widehat{G}_{1}(\ux)\widehat{\Phi}_{0}^{-}(\ux) + \widehat{G}^{*}_{0}(\ux)\widehat{\Phi}_{1}^{-}(\ux))= \widehat{g}_{1}(\ux),
\end{array}
\end{equation*}
which is equivalent to saying that
\begin{equation}\label{RBVP}
\begin{array}{ccccc}
\begin{pmatrix}
\widehat{\Phi}_{0}^{+}(\ux)\\ 
\widehat{\Phi}_{1}^{+}(\ux)
\end{pmatrix} & - & \begin{pmatrix}
\widehat{G}_{0}(\ux) & -\widehat{G}^{*}_{1}(\ux) \\ 
\widehat{G}_{1}(\ux) & \widehat{G}^{*}_{0}(\ux)
\end{pmatrix}\begin{pmatrix}
\widehat{\Phi}_{0}^{-}(\ux)\\ 
\widehat{\Phi}_{1}^{-}(\ux)
\end{pmatrix}  & = & \begin{pmatrix}
\widehat{g}_{0}(\ux)\\ 
\widehat{g}_{1}(\ux)
\end{pmatrix}.
\end{array} 
\end{equation}
Rewritten the Dirac operator $\partial_{n}$ in the form
\begin{equation*}
\partial_{n} := \sum_{i = 1}^{n}e_{i}\dfrac{\partial}{\partial x^{i}} = e_{1}(\dfrac{\partial}{\partial x^{1}} - \sum_{i = 2}^{n}e_{1}e_{i}\dfrac{\partial}{\partial x^{i}}) =: e_{1}\overline{D_{n - 1}^{'}}
\end{equation*}
\begin{equation*}
\partial_{n} := \sum_{i = 1}^{n}e_{i}\dfrac{\partial}{\partial x^{i}} = (\dfrac{\partial}{\partial x^{1}} +  \sum_{i = 2}^{n}e_{1}e_{i}\dfrac{\partial}{\partial x^{i}})e_{1} =: D_{n - 1}^{'} e_{1}.
\end{equation*}
we see
\begin{equation*}
\partial_{n}\Phi = e_{1}\overline{D_{n - 1}^{'}}(\Phi_{0}) - D_{n - 1}^{'}(\Phi_{1}) = - D_{n - 1}^{'}(\Phi_{1}) +  e_{1}\overline{D_{n - 1}^{'}}(\Phi_{0}).
\end{equation*}
Furthermore, the following are equivalent 
\begin{equation}\label{SistMon}
\begin{array}{ccc}
\partial_{n}\Phi = 0 & \Leftrightarrow  & \Big\lbrace \begin{array}{c}
D_{n - 1}^{'}(\Phi_{1}) = 0\\ 
\overline{D_{n - 1}^{'}}(\Phi_{0}) = 0.
\end{array} 
\end{array} 
\end{equation}
Besides this equivalence we have
\begin{equation*}
D_{n - 1} = \beta(D_{n - 1}^{'}) = \dfrac{\partial}{\partial x^{1}} +  \sum_{i = 1}^{n -1}e_{i}\dfrac{\partial}{\partial x^{i + 1}} = \partial_{1} + \partial_{n - 1}
\end{equation*}
\begin{equation*}
\overline{D_{n - 1}} = \beta(\overline{D_{n - 1}^{'}}) = \dfrac{\partial}{\partial x^{1}} - \sum_{i = 1}^{n -1}e_{i}\dfrac{\partial}{\partial x^{i + 1}} = \partial_{1} - \partial_{n - 1}
\end{equation*}
so that system ($\ref{SistMon}$) becomes
\begin{equation*}
\Big\lbrace \begin{array}{c}
D_{n - 1}(\widehat{\Phi}_{1}) = 0\\ 
\overline{D_{n - 1}}(\widehat{\Phi}_{0}) = 0.
\end{array} 
\end{equation*}
Summarizing, we have that $\widehat{\Phi}_{0}(x)$ and $\widehat{\Phi}_{1}(x)$ are antimonogenic and monogenic functions respectively. 
\begin{lemma}
A function is monogenic in the vectorial sense if and only if its even and odd part are, through isomorphism, antimonogenic and monogenic respectively in the paravectorial sense, and the following decomposition holds
\begin{equation*}
\Phi(\ux) = \beta^{-1}(\widehat{\Phi}_{0}(\ux)) + e_{1}\beta^{-1}(\widehat{\Phi}_{1}(\ux)),
\end{equation*}
where $\ux = \alpha^{-1}(x).$
\end{lemma}
This decomposition bears a striking similarity to that of analytic functions in complex analysis through two conjugate harmonic functions.

Since $\widehat{\Phi}_{0}(x)$ is left antimonogenic, then $\Upsilon_{0}(x):= \overline{\widehat{\Phi}_{0}(x)}$ is a right monogenic ones. For similarity we let $\Upsilon_{1} = \widehat{\Phi}_{1}$. Therefore problem (\ref{RBVP}) reduces to find 
\begin{equation*}
\begin{pmatrix}
\Upsilon_{0}(x)\\ 
\Upsilon_{1}(x)
\end{pmatrix},
\end{equation*}
such that on $\mathbb{R}^{n}\setminus L$
\begin{equation*}
\Big\lbrace \begin{array}{c}
(\Upsilon_{0})D_{n - 1} = 0\\ 
D_{n - 1}(\Upsilon_{1}) = 0,
\end{array} 
\end{equation*}
meanwhile on $L$ the boundary condition
\begin{equation}\label{CondComp2}
\begin{array}{ccccc}
\begin{pmatrix}
\overline{\Upsilon_{0}^{+}(\ux)}\\ 
\Upsilon_{1}^{+}(\ux)
\end{pmatrix} & - & \begin{pmatrix}
\widehat{G}_{0}(\ux) & -\widehat{G}^{*}_{1}(\ux) \\ 
\widehat{G}_{1}(\ux) & \widehat{G}^{*}_{0}(\ux)
\end{pmatrix}\begin{pmatrix}
\overline{\Upsilon_{0}^{-}(\ux)}\\ 
\Upsilon_{1}^{-}(\ux)
\end{pmatrix}  & = & \begin{pmatrix}
\widehat{g}_{0}(\ux)\\ 
\widehat{g}_{1}(\ux)
\end{pmatrix}
\end{array} 
\end{equation}
holds.

If this problem is solvable then so is ($\ref{ProbOrig}$) and explicit solution is given by
\begin{equation*}
\Phi(\ux) = \beta^{-1}\circ\overline{\Upsilon_{0}}\circ\alpha^{-1}(x) + e_{1}\beta^{-1}\circ\Upsilon_{1}\circ\alpha^{-1}(x) = 
\end{equation*}
\begin{equation}\label{SolCuater}
= \beta^{-1}(\overline{\Upsilon_{0}}(\ux)) + e_{1}\beta^{-1}(\Upsilon_{1}(\ux)),
\end{equation}
where
\begin{equation*}
\ux = \alpha^{-1}(x).
\end{equation*}
When we use the decomposition $\mathbb{R}_{0, n} = \mathbb{R}^{+}_{0, n}\bigoplus e_{n}\mathbb{R}^{+}_{0, n}$, analogous result can be obtained. We have
\begin{equation*}
\partial_{n} := \sum_{i = 1}^{n}e_{i}\dfrac{\partial}{\partial x_{i}} = (\dfrac{\partial}{\partial x_{n}} - \sum_{i = 1}^{n - 1}e_{i}e_{n}\dfrac{\partial}{\partial x_{i}})e_{n} =: \overline{D_{n - 1}^{''}}e_{n}
\end{equation*}
\begin{equation*}
\partial_{n} := \sum_{i = 1}^{n}e_{i}\dfrac{\partial}{\partial x_{i}} = e_{n}(\dfrac{\partial}{\partial x_{n}} +  \sum_{i = 1}^{n - 1}e_{i}e_{n}\dfrac{\partial}{\partial x_{i}}) =: e_{n}D_{n - 1}^{''}.
\end{equation*}
and thus
\begin{equation*}
\begin{array}{ccc}
\partial_{n}\Phi = 0 & \Leftrightarrow  & \Big\lbrace \begin{array}{c}
D_{n - 1}^{''}(\Phi_{0}) = 0\\ 
\overline{D_{n - 1}^{''}}(\Phi_{1}) = 0.
\end{array} 
\end{array} 
\end{equation*}
Being now
\begin{equation*}
\begin{array}{cccc}
\beta: & \mathbb{R}^{+}_{0, n}   & \rightarrow & \mathbb{R}_{0, n - 1} \\
       &     e_{i}e_{n}      & \rightarrow &  e_{i}
\end{array} 
\end{equation*}
\begin{equation*}
D_{n - 1} = \beta(\overline{D_{n - 1}^{''}}) = \dfrac{\partial}{\partial x_{n}} +  \sum_{i = 1}^{n -1}e_{i}\dfrac{\partial}{\partial x_{i}} = \partial_{1} + \partial_{n - 1}
\end{equation*}
\begin{equation*}
\overline{D_{n - 1}} = \beta(D_{n - 1}^{''}) = \dfrac{\partial}{\partial x_{n}} - \sum_{i = 1}^{n -1}e_{i}\dfrac{\partial}{\partial x_{i}} = \partial_{1} - \partial_{n - 1}
\end{equation*}
we have
\begin{equation*}
\Big\lbrace \begin{array}{c}
D_{n - 1}(\widehat{\Phi}_{0}) = 0\\ 
\overline{D_{n - 1}}(\widehat{\Phi}_{1}) = 0.
\end{array} 
\end{equation*}
\begin{lemma}
A function is monogenic in the vectorial sense if and only if its even and odd part are, through isomorphism, monogenic and antimonogenic respectively in the paravectorial sense, and the following decomposition holds
\begin{equation*}
\Phi(\ux) = \beta^{-1}(\widehat{\Phi}_{0}(\ux)) + e_{1}\beta^{-1}(\widehat{\Phi}_{1}(\ux)),
\end{equation*}
where $\ux = \alpha^{-1}(x) $.
\end{lemma}

\subsection{Cauchy type integral decomposition in vectorial Clifford analysis}
It has long been known that the RBVP theory in vectorial Clifford analysis is based on using the Cauchy type integral
\begin{equation*}
\Phi(\ux) = \int_{L}e_{n}(\underline{y} - \ux)\kappa(\underline{y})g(\underline{y})d\nu(\underline{y}).
\end{equation*}
In particular, for a smooth surface $L$, this integral, whose density $g$ satisfies a Holder condition, gives a unique solution to the simplest case of the RBVP (\ref{ProbOrig}), namely, the jump problem
\begin{equation*}
\Phi^{+}(\ux) - \Phi^{-}(\ux) = g(\ux),\,\ux\in L.
\end{equation*}
For what problem (\ref{CondComp2}) becomes
\begin{equation*}
\Big\lbrace \begin{array}{c}
\overline{\Upsilon_{0}^{+}(\ux)} - \overline{\Upsilon_{0}^{+}(\ux)} = \widehat{g}_{0}(\ux) \\
\Upsilon_{1}^{+}(\ux) -\Upsilon_{1}^{-}(\ux) = \widehat{g}_{1}(\ux),
\end{array} 
\end{equation*}
whose solution is
\begin{equation*}
\overline{\Upsilon_{0}(\ux)} = \overline{\int_{L}\overline{\widehat{g_{0}}(\uy)}\kappa(\uy)E_{n - 1}(\uy -\ux)d\nu(\uy)},
\end{equation*}
\begin{equation*}
\Upsilon_{1}(\ux) = \int_{L}E_{n - 1}(\uy - \ux)\kappa(\uy)\widehat{g_{1}}(\uy)d\nu(\uy).
\end{equation*}
Therefore by (\ref{SolCuater}) we obtain 
\begin{lemma}
The even part of the Cauchy type integral in vector Clifford analysis is the conjugate of a right Cauchy type integral in paravector Clifford analysis and its odd part is a left Cauchy type integral in paravector Clifford analysis both trough isomorphism.   
\end{lemma} 

\section{RBVP for monogenic functions in lower dimensional non-commutative Clifford algebras}
The lower dimensional non-commutative Clifford analysis focuses on functions $f: \mathbb{R}^{2} \rightarrow \mathbb{R}_{0, 2}$. Using
\begin{equation*}
\begin{array}{cccc}
\alpha: & \mathbb{C}   & \rightarrow & \mathbb{R}_{0,2}^{(1)}    \\ 
& x_{1} + x_{2}i & \rightarrow  & x_{1}e_{1} + x_{2}e_{2}   
\end{array} 
\end{equation*}
and
\begin{equation*}
\begin{array}{cccc}
\beta: & \mathbb{R}^{+}_{0, n}   & \rightarrow & \mathbb{C} \\
       &     e_{1}e_{2}      & \rightarrow &  i
\end{array} 
\end{equation*}
we can identify the correspondent paravector calculus with standard complex analysis. In fact
\begin{equation*}
\begin{array}{cccc}
\widehat{f} : & \mathbb{C}   & \rightarrow & \mathbb{C},   \\ 
 &      z := \alpha(\ux)     & \rightarrow &  \beta(f(\alpha(\ux))).
\end{array} 
\end{equation*}
Now, we can represent the Cauchy-Riemann operator and its conjugate as
\begin{equation*}
\partial_{\overline{z}} = \frac{1}{2}\beta(D_{1}^{'}) = \dfrac{1}{2}(\dfrac{\partial}{\partial x_{1}} + i\dfrac{\partial}{\partial x_{2}})
\end{equation*}
\begin{equation*}
\partial_{z} = \frac{1}{2}\beta(\overline{D_{1}^{'}}) = \dfrac{1}{2}(\dfrac{\partial}{\partial x_{1}} - i\dfrac{\partial}{\partial x_{2}}).
\end{equation*}

As a matter of fact,  $\widehat{G}^{*}_{j}(\ux) = \overline{\widehat{G}_{j}(\ux)}$, $j = 0, 1$. Then (\ref{CondComp2}) becomes 
\begin{equation}\label{CondComp3}
\begin{array}{ccccc}
\begin{pmatrix}
\overline{\Upsilon_{0}^{+}(\ux)}\\ 
\Upsilon_{1}^{+}(\ux)
\end{pmatrix} & - & \begin{pmatrix}
\widehat{G}_{0}(\ux) & -\overline{\widehat{G}_{1}(\ux)} \\ 
\widehat{G}_{1}(\ux) & \overline{\widehat{G}_{0}(\ux)}
\end{pmatrix}\begin{pmatrix}
\overline{\Upsilon_{0}^{-}(\ux)}\\ 
\Upsilon_{1}^{-}(\ux)
\end{pmatrix}  & = & \begin{pmatrix}
\widehat{g}_{0}(\ux)\\ 
\widehat{g}_{1}(\ux)
\end{pmatrix}.
\end{array} 
\end{equation}

\subsection{Case of null odd part}
We have $G_{1} \equiv 0$ and hence $G = G_{0} \neq 0$. Consequently, ($\ref{CondComp3}$) becomes 
\begin{equation*}
\Big\lbrace \begin{array}{c}
\overline{\Upsilon_{0}^{+}(\ux)} - \widehat{G}_{0}(\ux)\overline{\Upsilon_{0}^{+}(\ux)} = \widehat{g}_{0}(\ux) \\
\Upsilon_{1}^{+}(\ux) -\overline{\widehat{G}_{0}(\ux)}\Upsilon_{1}^{-}(\ux) = \widehat{g}_{1}(\ux),
\end{array} 
\end{equation*}
which will lead to
\begin{equation*}
\Big\lbrace \begin{array}{c}
\Upsilon_{0}^{+}(\ux) - \overline{\widehat{G}_{0}(\ux)}\Upsilon_{0}^{+}(\ux) = \overline{\widehat{g}_{0}(\ux)} \\
\Upsilon_{1}^{+}(\ux) -\overline{\widehat{G}_{0}(\ux)}\Upsilon_{1}^{-}(\ux) = \widehat{g}_{1}(\ux).
\end{array} 
\end{equation*}
These are two independent RBVPs in complex analysis with the same coefficient. We have $\overline{\widehat{G}_{0}(\ux)} \neq 0$ since $G(\ux) \neq 0$.

We define $Ind(G):= Ind(\widehat{G}_{0})= \aleph$, which yields $Ind(\overline{\widehat{G}_{0}}) = -\aleph$. Following the standard techniques of the RBVP theory we have
\begin{equation*}
\begin{array}{cc}
X^{+}(z) = e^{\varGamma(z)},	& X^{-}(z) = z^{\aleph}e^{\varGamma(z)},
\end{array} 
\end{equation*}
where
\begin{equation*}
\varGamma(z) =\dfrac{1}{2\pi i}\int\limits_{L}\dfrac{\ln[\tau^{\aleph}\overline{\widehat{G}_{0}(\tau)}]}{\tau - z}d\tau,
\end{equation*}
and
\begin{equation*}
\Psi_{0}(z) =\dfrac{1}{2\pi i}\int\limits_{L}\dfrac{\overline{\widehat{g}_{0}(\tau)}}{X^{+}(\tau)}\dfrac{d\tau}{\tau - z}.
\end{equation*}

\begin{equation*}
\Psi_{1}(z) =\dfrac{1}{2\pi i}\int\limits_{L}\dfrac{\widehat{g}_{1}(\tau)}{X^{+}(\tau)}\dfrac{d\tau}{\tau - z}.
\end{equation*}

\begin{theorem}
If $\aleph \leq 1$, the solution $\Phi(\ux)$ is obtained by (\ref{SolCuater}) where 
\begin{equation*}
\Upsilon_{0}(z) = X(z)[\Psi_{0}(z) + P^{0}_{-\aleph}(z)],	
\end{equation*} 	
\begin{equation*}
\Upsilon_{1}(z) = X(z)[\Psi_{1}(z) + P^{1}_{-\aleph}(z)].	
\end{equation*} 	
Here $P^{0}_{-\aleph}(z), P^{1}_{-\aleph}(z)$ are two polynomials of degree $-\aleph$. For $\aleph = 1$ we put $P^{0}_{-\aleph}(z) \equiv 0, P^{1}_{-\aleph}(z) \equiv 0$.

If $\aleph > 1 $, when the following $2(\aleph - 1)$ solubility conditions
\begin{equation*}
\begin{array}{cc}
\int\limits_{L}\dfrac{\overline{\widehat{g}_{0}(\tau)}}{X^{+}(\tau)}\tau^{k-1}d\tau = 0, & k = 1, 2, \ldots, \aleph - 1
\end{array} 
\end{equation*}
\begin{equation*}
\begin{array}{cc}
\int_{L}\dfrac{\widehat{g}_{1}(\tau)}{X^{+}(\tau)}\tau^{k-1}d\tau = 0, & k = 1, 2, \ldots, \aleph - 1
\end{array} 
\end{equation*}
are fulfilled, thus (\ref{SolCuater}) is the solution, where $P^{0}_{-\aleph}(z) \equiv 0, P^{1}_{-\aleph}(z) \equiv 0$.
\end{theorem}
Under condition $\Phi^{-}(\infty) = 0$, our theorem get a more symmetrical form
\begin{theorem}
Under the condition $\Phi^{-}(\infty) = 0$, if $\aleph \leq 0$, the solution $\Phi(\ux)$ is given by (\ref{SolCuater}), where 
\begin{equation*}
\Upsilon_{0}(z) = X(z)[\Psi_{0}(z) + P^{0}_{-\aleph - 1}(z)],	
\end{equation*} 	
\begin{equation*}
\Upsilon_{1}(z) = X(z)[\Psi_{1}(z) + P^{1}_{-\aleph - 1}(z)].	
\end{equation*} 	
Here $P^{0}_{-\aleph - 1}(z), P^{1}_{-\aleph - 1}(z)$ are two polynomials of degree $-\aleph - 1$. For $\aleph = 0$ we put $P^{0}_{-\aleph - 1}(z) \equiv 0, P^{1}_{-\aleph - 1}(z) \equiv 0$, and the solution depends of $-4\aleph$ real constants.

If $\aleph > 0 $, when the following $2\aleph$ solubility conditions 
\begin{equation*}
\begin{array}{cc}
\int_{L}\dfrac{\overline{\widehat{g}_{0}(\tau)}}{X^{+}(\tau)}\tau^{k-1}d\tau = 0, & k = 1, 2, \ldots, \aleph 
\end{array} 
\end{equation*}
\begin{equation*}
\begin{array}{cc}
\int_{L}\dfrac{\widehat{g}_{1}(\tau)}{X^{+}(\tau)}\tau^{k-1}d\tau = 0, & k = 1, 2, \ldots, \aleph
\end{array} 
\end{equation*}
are fulfilled, thus (\ref{SolCuater}) is the solution, where $P^{0}_{-\aleph - 1}(z) \equiv 0, P^{1}_{-\aleph - 1}(z) \equiv 0$.
\end{theorem}

\subsubsection*{Example 1}
Find a function $\Phi$, vanishing at infinity, holomorphic in $\R^2\setminus L$ by the boundary condition
\begin{equation*}
\Phi^{+}(\ux) - (-e_{1}\ux)[(-e_{1}\ux)^{2} - 1]^{-1}\Phi^{-}(\ux) = [-e_{1}\ux - 1]^{-1} + e_{1}[-\ux e_{1} + 1]^{-1},\,\ux\in L,
\end{equation*}
where $L$ is an arbitrary smooth curve assuming additional conditions to be divided in four cases.

In this example, the problem reduces to the complex Riemann problems:
\begin{equation*}
\Big\lbrace \begin{array}{c}
\Upsilon_{0}^{+}(\ux) - \dfrac{t}{t^{2} - 1}\Upsilon_{0}^{+}(\ux) =\dfrac{1}{t - 1}\\
\Upsilon_{1}^{+}(\ux) - \dfrac{t}{t^{2} - 1}\Upsilon_{1}^{-}(\ux) = \dfrac{1}{t + 1}
\end{array} 
\end{equation*}

\noindent
\textbf{Case a)} $L$ contains inside the point $z = 0$ and $z = 1, z = -1$ are outside. Let $\aleph = Ind(G) = -Ind(\overline{\widehat{G}_{0}(\ux)}) = -1$, so that
\begin{equation*}
\begin{array}{cc}
X^{+}(z) = \dfrac{1}{z^{2} - 1},	& X^{-}(z) = \dfrac{1}{z},\\
\Psi^{+}_{0}(z) = z + 1, & \Psi^{-}_{0}(z) = 0, \\
\Psi^{+}_{1}(z) = z - 1, & \Psi^{-}_{1}(z) = 0. \\
\end{array}
\end{equation*}
We have
\begin{equation*}
\begin{array}{cc}
\Upsilon_{0}^{+}(z) = \dfrac{1}{z^{2} - 1}(z + 1 + c_{1}), & \Upsilon_{0}^{-}(z) = \dfrac{c_{1}}{z},\\
\Upsilon_{1}^{+}(z) = \dfrac{1}{z^{2} - 1}(z - 1 + c_{2}), & \Upsilon_{1}^{-}(z) = \dfrac{c_{2}}{z},
\end{array}
\end{equation*}
and thus
\begin{equation*}
\Phi^{+}(\ux) = [-e_{1}\ux - 1]^{-1} + (c^{0}_{1} + e_{1}e_{2}c^{1}_{1})[(-e_{1}\ux)^{2} - 1]^{-1} + e_{1}[-\ux e_{1} + 1]^{-1} + 
\end{equation*}
\begin{equation*}
 + (e_{1}c^{0}_{2} -e_{2}c^{1}_{2})[(-\ux e_{1})^{2} - 1]^{-1},
\end{equation*}
\begin{equation*}
\Phi^{-}(z) = (c^{0}_{1} + e_{1}e_{2}c^{1}_{1})[-e_{1}\ux]^{-1} + (e_{1}c^{0}_{2} - e_{2}c^{1}_{2})[-\ux e_{1}]^{-1},
\end{equation*}

\noindent
\textbf{Case b)} $L$ contains inside the points $z = 0, z = 1$ and $z = -1$ lies outside. Then $\aleph = Ind(G) = -Ind(\overline{\widehat{G}_{0}(\ux)}) = 0$, so that
\begin{equation*}
\begin{array}{cc}
X^{+}(z) = \dfrac{1}{z + 1},	& X^{-}(z) = \dfrac{z - 1}{z},\\
\Psi^{+}_{0}(z) = 1, & \Psi^{-}_{0}(z) = -\dfrac{2}{z - 1}, \\
\Psi^{+}_{1}(z) = 1, & \Psi^{-}_{1}(z) = 0. \\
\end{array}
\end{equation*}
We have
\begin{equation*}
\begin{array}{cc}
\Upsilon_{0}^{+}(z) = \dfrac{1}{z + 1}, & \Upsilon_{0}^{-}(z) = -\dfrac{2}{z},\\
\Upsilon_{1}^{+}(z) = \dfrac{1}{z + 1}, & \Upsilon_{0}^{-}(z) = 0,
\end{array}
\end{equation*}
and 
\begin{equation*}
\begin{array}{cc}
\Phi^{+}(\ux) = [-e_{1}\ux + 1]^{-1} + e_{1}[-\ux e_{1} + 1]^{-1}, & \Phi^{-}(z) = -2[-e_{1}\ux]^{-1}.
\end{array}
\end{equation*}

\noindent
\textbf{Case c)} $L$ contains inside the point $z = -1$ and $z = 0, z = 1$ are outside. Thus $\aleph = Ind(G) = -Ind(\overline{\widehat{G}_{0}(\ux)}) = 1$. Then the solubility conditions 
\begin{equation*}
\int_{L}\dfrac{1}{\tau}d\tau = 0, 
\end{equation*}
\begin{equation*}
\int_{L}\dfrac{\tau - 1}{\tau(\tau + 1)}d\tau = \int_{L}\dfrac{1}{\tau + 1}d\tau - \int_{L}\dfrac{\frac{1}{\tau}}{\tau + 1}d\tau = 4\pi i, 
\end{equation*}
must be satisfied. So the problem has no solution. Notice that, in this case, the first complex problem has a solution, but it is not enough to solve the original problem.
 
\noindent
\textbf{Case d)} $L$ contains inside the point $z = -1$ and $z = 0, z = 1$ are outside. Let $\aleph = Ind(G) = -Ind(\overline{\widehat{G}_{0}(\ux)}) = 1$,  but in this case function $g(\ux)$ should be taking quite different
\begin{equation*}
g(\ux) = [-e_{1}\ux - 1]^{-1} + e_{1}[-\ux e_{1}]^{-1}.
\end{equation*}
Checking the solubility conditions we have
\begin{equation*}
\int_{L}\dfrac{1}{\tau}d\tau = 0, 
\end{equation*}
\begin{equation*}
\int_{L}\dfrac{\tau - 1}{\tau^{2}}d\tau = \int_{L}\dfrac{1}{\tau}d\tau - \int_{L}\dfrac{1}{\tau^{2}}d\tau = 0. 
\end{equation*}
Therefore, there exists solution
\begin{equation*}
\begin{array}{cc}
X^{+}(z) = \dfrac{z}{z - 1},	& X^{-}(z) = z + 1,\\
\Psi^{+}_{0}(z) = \dfrac{1}{z}, & \Psi^{-}_{0}(z) = 0, \\
\Psi^{+}_{1}(z) = \dfrac{z - 1}{z^{2}}, & \Psi^{-}_{1}(z) = 0. \\
\end{array}
\end{equation*}
We have
\begin{equation*}
\begin{array}{cc}
\Upsilon_{0}^{+}(z) = \dfrac{1}{z - 1}, & \Upsilon_{0}^{-}(z) = 0,\\
\Upsilon_{1}^{+}(z) = \dfrac{1}{z}, & \Upsilon_{0}^{-}(z) = 0.
\end{array}
\end{equation*}
Thus, the only solution is
\begin{equation*}
\begin{array}{cc}
\Phi^{+} = [-e_{1}\ux - 1]^{-1} + e_{1}[-\ux e_{1}]^{-1}, & \Phi^{-} = 0.
\end{array}
\end{equation*}

\subsection{Case of constant coefficients}
For this case, the theory of conformal mappings will be used. Let $L$ denote a simple closed and smooth contour with a tangent that forms a certain angle with a constant direction that satisfies a H\"{o}lder condition. This idea has previously been used in \cite{Chibrikova-Rogozhin, Mu}. 

\subsubsection{Reduction to a circle case}
Let us denote by $\omega = \chi^{+}(z)$ ($\omega = \chi^{-}(z)$) a conformal mapping from $\Omega_{+}$ ($\Omega_{-}$) to the inside (outside) of the unit circle $C$. We shall write $z = \varphi^{+}(\omega)$ ($z = \varphi^{-}(\omega)$) the respective inverses. From the theory of conformal mappings it is known that under the adopted conditions referred to the contour $L$ not only the functions $\chi^{+}(z), \chi^{-}(z),\varphi^{+}(\omega), \varphi^{-}(\omega)$, but also its first derivatives are continuously prolonged over $L$ and $C$ respectively and satisfy a H\"older condition.

Introducing the new functions
\begin{equation*}
\begin{array}{ccc}
\Psi^{+}_{j}(\omega) = \Upsilon^{+}_{j}[\varphi^{+}(\omega)],	& \Psi^{-}_{j}(\omega) = \Upsilon^{-}_{j}[\varphi^{-}(\omega)], & j = 0, 1,
\end{array} 
\end{equation*}
the boundary condition takes the form
\begin{equation}\label{ProbCir}
\begin{array}{ccccc}
\begin{pmatrix}
\overline{\Psi^{+}_{0}(\zeta)}\\ 
\Psi^{+}_{1}(\zeta)
\end{pmatrix} & - & \begin{pmatrix}
\widetilde{G}_{0}(\zeta) & -\overline{\widetilde{G}_{1}(\zeta)} \\ 
\widetilde{G}_{1}(\zeta) & \overline{\widetilde{G}_{0}(\zeta)}
\end{pmatrix}\begin{pmatrix}
\overline{\Psi_{0}^{-}(\zeta)}\\ 
\Psi_{1}^{-}(\zeta)
\end{pmatrix}  & = & \begin{pmatrix}
\widetilde{g}_{0}(\zeta)\\ 
\widetilde{g}_{1}(\zeta)
\end{pmatrix}.
\end{array} 
\end{equation}
where $\widetilde{g}_{j}(\zeta) = \widehat{g}_{j}(\varphi^{-}(\zeta))$, $\widetilde{G}_{j}(\zeta) = \widehat{G}_{j}(\varphi^{-}(\zeta))$, $j = 0, 1$.\\
The functions $\widetilde{G}_{j}, \widetilde{g}_{j}$ defined on $C$ satisfy a H\"older condition, when it is fulfilled by $\widehat{G}_{j}, \widehat{g}_{j}$. So the problem (\ref{CondComp2}) becomes at (\ref{ProbCir}) considered on the unit circle with center at the origin.

\subsubsection{Solution over the unit circle}
Let $L$ be a unit circle with center at the origin. Doing in (\ref{CondComp2}) the change of variables
\begin{equation*}
\Lambda^{+}(z) = \overline{\Upsilon_{0}^{-}(\frac{1}{\overline{z}})},
\end{equation*}
\begin{equation*}
\Lambda^{-}(z) = \overline{\Upsilon_{0}^{+}(\frac{1}{\overline{z}})},
\end{equation*}
note that $\overline{\Upsilon_{0}^{-}(\infty)} = \Lambda^{+}(0)$. If we have the condition $\Phi(\infty) = 0$, then for (\ref{SolCuater}) we have that $\Lambda^{+}(0) = 0$,  so over $L$ we have that
\begin{equation*}
\Lambda^{+}(\ux) = \overline{\Upsilon_{0}^{-}(\frac{1}{\overline{t}})} = \overline{\Upsilon_{0}^{-}(\ux)},
\end{equation*}
\begin{equation*}
\Lambda^{-}(\ux) = \overline{\Upsilon_{0}^{+}(\frac{1}{\overline{t}})} = \overline{\Upsilon_{0}^{+}(\ux)}.
\end{equation*}
next the system (\ref{CondComp2}) becomes 
\begin{equation*}
\begin{array}{ccccc}
\begin{pmatrix}
\Lambda^{-}(\ux)\\ 
\Upsilon_{1}^{+}(\ux)
\end{pmatrix} & - & \begin{pmatrix}
\widehat{G}_{0}(\ux) & -\overline{\widehat{G}_{1}(\ux)} \\ 
\widehat{G}_{1}(\ux) & \overline{\widehat{G}_{0}(\ux)}
\end{pmatrix}\begin{pmatrix}
\Lambda^{+}(\ux)\\ 
\Upsilon_{1}^{-}(\ux)
\end{pmatrix}  & = & \begin{pmatrix}
\widehat{g}_{0}(\ux)\\ 
\widehat{g}_{1}(\ux)
\end{pmatrix}
\end{array} 
\end{equation*}
where the function $\Lambda$ is holomorphic on $\C\setminus L$.

Now if we consider the case of constant coefficients we have that $\widehat{G}_{0}(\ux) \equiv a$ and $\widehat{G}_{1}(\ux) \equiv b$ where the constants $a$, $b$ $\in \mathbb{C}$ and then we obtain that
\begin{equation*}
\begin{array}{ccccc}
\begin{pmatrix}
\Lambda^{-}(\ux)\\ 
\Upsilon_{1}^{+}(\ux)
\end{pmatrix} & - & \begin{pmatrix}
a & -\overline{b} \\ 
b & \overline{a}
\end{pmatrix}\begin{pmatrix}
\Lambda^{+}(\ux)\\ 
\Upsilon_{1}^{-}(\ux)
\end{pmatrix}  & = & \begin{pmatrix}
\widehat{g}_{0}(\ux)\\ 
\widehat{g}_{1}(\ux)
\end{pmatrix}
\end{array} 
\end{equation*}
that is
\begin{equation}\label{CondConsCom1}
\Lambda^{-}(\ux) - a\Lambda^{+}(\ux) + \overline{b}\Upsilon_{1}^{-}(\ux) = \widehat{g}_{0}(\ux),
\end{equation}
\begin{equation}\label{CondConsCom2}
\Upsilon_{1}^{+}(\ux) - b\Lambda^{+}(\ux) - \overline{a}\Upsilon_{1}^{-}(\ux) = \widehat{g}_{1}(\ux).
\end{equation}
Since in (\ref{ProbOrig}) we have $G(\ux) \neq 0$ then at least $a \neq$ or $b \neq 0$. We will consider three cases.

\noindent 
\textbf{Case a)}We will analyze first the case in which $a \neq 0$ and $b \neq 0$. Here we have for (\ref{CondConsCom1})
\begin{equation*}
\Lambda^{+}(\ux) -\frac{1}{a}\Lambda^{-}(\ux) - \frac{\overline{b}}{a}\Upsilon_{1}^{-}(\ux) = -\frac{1}{a}\widehat{g}_{0}(\ux),
\end{equation*}
and put 
\begin{equation*}
\begin{array}{cc}
\Psi_{0}^{+}(z) = \Lambda^{+}(z)      & \Psi_{0}^{-}(z) = \frac{1}{a}\Lambda^{-}(z) + \frac{\overline{b}}{a}\Upsilon_{1}^{-}(z)
\end{array}
\end{equation*}
the problem becomes in
\begin{equation*}
\Psi_{0}^{+}(\ux) -\Psi_{0}^{-}(\ux) = -\frac{1}{a}\widehat{g}_{0}(\ux),
\end{equation*}
whit $\Psi_{0}^{+}(0) = 0$. The solution to this problem is
\begin{equation*}
\Psi_{0}(z) = \frac{1}{2\pi i}\int_{L}\dfrac{-\frac{1}{a}\widehat{g}_{0}(\tau)}{\tau - z} - \frac{1}{2\pi i}\int_{L}\dfrac{-\frac{1}{a}\widehat{g}_{0}(\tau)}{\tau}.
\end{equation*}
Substituting in (\ref{CondConsCom2})
\begin{equation*}
\Upsilon_{1}^{+}(\ux)  - \overline{a}\Upsilon_{1}^{-}(\ux) = \widehat{g}_{1}(\ux) + b\Psi_{0}^{+}(\ux).
\end{equation*}
and put 
\begin{equation*}
\begin{array}{cc}
\Psi_{1}^{+}(z) = \Upsilon_{1}^{+}(z)      & \Psi_{1}^{-}(z) = \overline{a}\Upsilon_{1}^{-}(z)
\end{array}
\end{equation*}
the problem becomes in
\begin{equation*}
\Psi_{1}^{+}(\ux) -\Psi_{1}^{-}(\ux) = \widehat{g}_{1}(\ux) + b\Psi_{0}^{+}(\ux).
\end{equation*}
whit $\Psi_{1}^{-}(\infty) = 0$. The solution to this problem is
\begin{equation*}
\Psi_{1}(z) = \frac{1}{2\pi i}\int_{L}\dfrac{\widehat{g}_{1}(\tau) + b\Psi_{0}^{+}(\tau)}{\tau - z}.
\end{equation*}
Therefore we obtain
\begin{equation*}
\begin{array}{cc}
\Lambda^{+}(z) = \Psi_{0}^{+}(z) & \Lambda^{-}(z) = a\Psi_{0}^{-}(z) - \frac{\overline{b}}{\overline{a}}\Psi_{1}^{-}(z) \\ 
\Upsilon_{1}^{+}(z) = \Psi_{1}^{+}(z) & \Upsilon_{1}^{-}(z) = \frac{1}{\overline{a}}\Psi_{1}^{-}(z). 
\end{array} 
\end{equation*}

\subsubsection*{Example 2}
Find a function $\Phi$, that vanishes at infinity, holomorphic in $\R^2\setminus L$ and satisfying:
\begin{equation*}
\Phi^{+}(\ux) - (1 - e_{1})\Phi^{-}(\ux) = [-\ux e_{1}]^{-1} + e_{1}[-\ux e_{1} - 2]^{-1},\,\ux\in L.
\end{equation*}
where $L$ is the unit circle with center in the origin.\\
In this case we have
\begin{equation*}
\Lambda^{-}(\ux) - \Lambda^{+}(\ux) + \Upsilon_{1}^{-}(\ux) = \dfrac{1}{t},
\end{equation*}
\begin{equation*}
\Upsilon_{1}^{+}(\ux) - \Lambda^{+}(\ux) - \Upsilon_{1}^{-}(\ux) = \dfrac{1}{t - 2}.
\end{equation*}
and thus
\begin{equation*}
\begin{array}{cc}
\Psi_{0}^{+}(z) \equiv 0 & \Psi_{0}^{-}(z) = \dfrac{1}{z} \\ 
\Psi_{1}^{+}(z) = \dfrac{1}{z - 2}  & \Psi_{1}^{-}(z) \equiv 0. 
\end{array} 
\end{equation*}
then we have
\begin{equation*}
\begin{array}{cc}
\Lambda^{+}(z) \equiv 0 & \Lambda^{-}(z) \equiv \dfrac{1}{z} \\ 
\Upsilon_{1}^{+}(z) = \dfrac{1}{z - 2} & \Upsilon_{1}^{-}(z) \equiv  0. 
\end{array} 
\end{equation*}
next
\begin{equation*}
\overline{\Upsilon_{0}^{-}(z)} = \Lambda^{+}(\frac{1}{\overline{z}}) \equiv 0,
\end{equation*}
\begin{equation*}
\overline{\Upsilon_{0}^{+}(z)} = \Lambda^{-}(\frac{1}{\overline{z}}) = \overline{z} ,
\end{equation*}
therefore the solution is
\begin{equation*}
\begin{array}{cc}
\Phi^{+}(\ux) = (-e_{1}\ux) + e_{1}[-\ux e_{1} - 2]^{-1}  & \Phi^{-}(\ux) = 0.
\end{array} 
\end{equation*}

\noindent
\textbf{Case b)} $a \neq 0$ and $b = 0$ this can be handled in much the same way. Indeed we have 
\begin{equation*}
\begin{array}{cc}
\Lambda^{+}(z) = \Psi_{0}^{+}(z) & \Lambda^{-}(z) = a\Psi_{0}^{-}(z) \\ 
\Upsilon_{1}^{+}(z) = \Psi_{1}^{+}(z) & \Upsilon_{1}^{-}(z) = \frac{1}{\overline{a}}\Psi_{1}^{-}(z). 
\end{array} 
\end{equation*}
where
\begin{equation*}
\Psi_{0}(z) = \frac{1}{2\pi i}\int_{L}\dfrac{-\frac{1}{a}\widehat{g}_{0}(\tau)}{\tau - z} - \frac{1}{2\pi i}\int_{L}\dfrac{-\frac{1}{a}\widehat{g}_{0}(\tau)}{\tau}.
\end{equation*}
\begin{equation*}
\Psi_{1}(z) = \frac{1}{2\pi i}\int_{L}\dfrac{\widehat{g}_{1}(\tau)}{\tau - z}.
\end{equation*}

\noindent
\textbf{Case c)} $a = 0$ and $b \neq 0$ here we have
\begin{equation*}
\Lambda^{-}(\ux) + \overline{b}\Upsilon_{1}^{-}(\ux) = \widehat{g}_{0}(\ux),
\end{equation*}
\begin{equation*}
\Upsilon_{1}^{+}(\ux) - b\Lambda^{+}(\ux) = \widehat{g}_{1}(\ux).
\end{equation*}
A necessary and sufficient condition for this problem to be solved is that $\widehat{g}_{0}(\ux)$ admits a holomorphic extension to $\Omega^{+}$ and $\widehat{g}_{0}(\ux)$. If these are satisfied we can choose one of this functions for example $\Lambda(z)$ being holomorphic in $\mathbb{C}\backslash L$ and $\Lambda(0) = 0$. In particular we can choose $\Lambda^{+}(z) = z^{n}$ and $\Lambda^{-}(z) = \displaystyle\frac{1}{z^{n}}$ and then we have  
\begin{equation*}
\Upsilon_{1}^{+}(z) = \frac{1}{2\pi i}\int_{L}\dfrac{\widehat{g}_{1}(\tau) + b\tau^{n}}{\tau - z},
\end{equation*}
\begin{equation*}
\Upsilon_{1}^{-}(z) = -\frac{1}{2\overline{b}\pi i}\int_{L}\dfrac{\widehat{g}_{0}(\tau) - \frac{1}{\tau^{n}}}{\tau - z}.
\end{equation*}
Next in the tree cases doing
\begin{equation*}
\overline{\Upsilon_{0}^{-}(z)} = \Lambda^{+}(\frac{1}{\overline{z}}),
\end{equation*}
\begin{equation*}
\overline{\Upsilon_{0}^{+}(z)} = \Lambda^{-}(\frac{1}{\overline{z}}),
\end{equation*}
and applying the reverse mappings through (\ref{SolCuater}) we have the solution. We have obtained the following theorem.
\begin{theorem}
If the even part of the Riemann boundary value problem with constants coefficients is no null then the problem have unique solution else if not the following solubility conditions must be satisfied
\begin{equation*}
\frac{1}{2}\widehat{g}_{0}(\ux) + \frac{1}{2\pi i}\int_{L}\dfrac{\widehat{g}_{0}(\tau)}{\tau - t} = 0
\end{equation*}
\begin{equation*}
-\frac{1}{2}\widehat{g}_{1}(\ux) + \frac{1}{2\pi i}\int_{L}\dfrac{\widehat{g}_{1}(\tau)}{\tau - t} = 0
\end{equation*}
if both are satisfied then the problem has an infinite number of linearly independent solutions, that vanishes at infinity.
\end{theorem}

\subsubsection*{Example 3}
Find a function $\Phi$, that vanishes at infinity, holomorphic in $\R^2\setminus L$ and satisfying:
\begin{equation*}
\Phi^{+}(\ux) - e_{1}\Phi^{-}(\ux) = 0,\,\ux\in L.
\end{equation*}
where $L$ is the unit circle with center in the origin.\\
In this case we have
\begin{equation*}
\Lambda^{-}(\ux) + \Upsilon_{1}^{-}(\ux) = 0,
\end{equation*}
\begin{equation*}
\Upsilon_{1}^{+}(\ux) - \Lambda^{+}(\ux) = 0.
\end{equation*}
where $\Lambda^{+}(0) = 0$ and $\Upsilon_{1}^{-}(\infty) = 0$.\\
Obviously the solubility conditions are satisfied, then we can choose $\Lambda^{+}(z) = 0$ and $\Lambda^{-}(z) = \displaystyle\frac{1}{z^{m}}, m\in {\mathbb N}$ and then we have 
\begin{equation*}
\begin{array}{cc}
\overline{\Upsilon_{0}^{+}(z)} = \overline{z}^{m}, & \overline{\Upsilon_{0}^{-}(z)} = 0, \\ 
\Upsilon_{1}^{+}(z) = 0, & \Upsilon_{1}^{-}(z) = -\displaystyle\frac{1}{z^{m}}, m\in {\mathbb N}. 
\end{array} 
\end{equation*}
Therefore we obtain
\begin{equation*}
\begin{array}{cc}
\Phi^{+}(\ux) = (-e_{1}\ux)^{m}, & \Phi^{-}(\ux) = - e_{1}[(-\ux e_{1})^{m}]^{-1}, m\in {\mathbb N}.
\end{array} 
\end{equation*}
By decreasing induction on $m$, we can verify that these functions satisfy the conditions of the problem.

The significance of this example is captured by the following corollary, which is consistent with the earlier results on the Fredholmness of the left linear Riemann operator reported in \cite{Br,ShV}.
\begin{corollary}
The homogeneous Riemann boundary value problem with constants coefficients can has an infinite number of linearly independent solutions, that vanishes at infinity. 
\end{corollary}

\section*{Acknowledgements}
C. D. Tamayo Castro gratefully acknowledges the financial support of the Postgraduate Study Fellowship of the Consejo Nacional de Ciencia y Tecnolog\'ia (CONACYT) (grant number 957110). J. Bory Reyes was partially supported by Instituto Polit\'ecnico Nacional in the framework of SIP programs (SIP20195662).

\end{document}